\documentclass[12pt]{amsart}
\usepackage{amsfonts}

\theoremstyle{definition}

\theoremstyle{remark}

\newcommand{\const}{\mathop{\rm const}\limits}

\newcommand{\Ent}{\mathop{\rm Ent}\limits}

\begin{document}

\begin{center}

{\bf BESOV'S TYPE EMBEDDING THEOREM FOR BILATERAL GRAND LEBESGUE SPACES} \par

\vspace{3mm}
{\bf E. Ostrovsky and L. Sirota}\\

\vspace{2mm}

{\it Department of Mathematics and Statistics, Bar-Ilan University,
59200, Ramat Gan.}\\
e \ - \ mails: galo@list.ru; \ sirota@zahav.net.il \\

\vspace{4mm}

 {\sc Abstract.} \\

\end{center}

{\it In this paper we obtain the non - asymptotic norm estimations of Besov's type
between  the norms of a functions in different Bilateral Grand  Lebesgue spaces
(BGLS).\par
We also give some examples to show the sharpness of these inequalities.} \\

\vspace{3mm}


\vspace{3mm}


\vspace{3mm}

\section{Introduction. Statement of problem. Notations. Auxiliary facts.}

\vspace{3mm}
{\bf A. Introduction.} \par
Let $ X $ be the set $ [0, 2 \pi] $ endowed with the normalized Lebesgue measure

$$
\mu(A) = (2 \pi)^{-1} \int_A dx, \ A \subset X.
$$
 Let also $ B $ be some rearrangement invariant (r.i.) space with the norm of a function
 $ f: X \to R $ denoted by $ ||f||B $ over the set $ X. $ We denote by
 $ T(h)f $ for arbitrary function $ f $ from the space $ B $ the shift operator
$$
T(h)f(x) = f( x+h),  h \in (-2\pi, 2 \pi),
$$
where all arithmetical operations over the arguments of a functions are understood
modulo $ 2 \pi; $ and we introduce the correspondent modulo of continuity

$$
\omega^{(B)}(f,\delta) \stackrel{def}{=} \sup_{h:|h|\le \delta}||T(h)f - f||B,  \
\delta \in X. \eqno(1)
$$

\vspace{2mm}

 We denote as usually the classical $ L(p), \ p \ge 1 $ Lebesgue norm

 $$
 |f|_p = \left( \int_{X} |f(x)|^p \ d\mu  \right)^{1/p};  \  f \in L(p) \
 \Leftrightarrow |f|_p < \infty.
 $$
We will denote for the spaces $ L(p) $

$$
\omega^{(L(p))}(f,\delta) = \omega(f,\delta)_p.
$$
 It follows from H\"older inequality that if $  1 \le p \le q, $ then
 $ |f|_p \le |f|_q,$ or simple $ L(q) \subset L(p).$ \par
  Inversely, if $ f \in L(p) $ and
$$
\omega^{(L(p))}(f,\delta) = \omega(f,\delta)_p \le C \ \delta^{\alpha},\eqno(2)
$$
where $C, \alpha = \const >0, \ \alpha p < 1, $ then for any $ q \in (p, p/(1-\alpha p)) $
$$
f \in L(q),  \eqno(3)
$$
see, e.g. \cite{Besov1}, chapter 4; \cite{Besov2}, chapter 2.\par
 Analogous assertions are true for the classical Besov's spaces \cite{Besov1}, \cite{Besov2},  \cite{Bennet1}, chapter 5, section 4, p. 332-346 etc;
 for the Sobolev's spaces \cite{Adams1}, \cite{Chuas1},  \cite{Edmunds1},
  \cite{Kufner1}, \cite{Maz'ja1}, \cite{Perez1} \cite{Stein1},
\cite{Talenty1}, \cite{Wannebo1}   etc.  \par

\vspace{2mm}

{\bf B. Statement of problem.}\par

\vspace{3mm}

{\it  Our aim is a generalization of the estimation (2), (3) on the so-called Bilateral
 Grand Lebesgue Spaces $ BGL = BGL(\psi) = G(\psi), $ i.e. when } $ f(\cdot)
 \in G(\psi) \ $  {\it and to show the precision of obtained estimations by means
 of the constructions of suitable examples. } \par

 \vspace{3mm}

 {\bf C. Another notations.}\par
 \vspace{3mm}
 Let $ f \in L(p), \ p \in (1,\infty) $  and $ n = 1,2,3,\ldots. $
 We denote as ordinary the error of the best trigonometric approximation in the norm
 $ |\cdot|_p $  of the function $ f(\cdot) $ by the trigonometric polynomials of
 degree $ \le n $ by $ E_n(f)_p: $

 $$
 E_n(f)_p = \inf_{\deg(U) \le n} |f - U|_p,
 $$
 where $ U = U(x)$ is arbitrary trigonometric polynomial of degree less or equal $ n: $
 $ U(\cdot) \in Q(n), $ where $ Q(n) $ is the {\it set} of all
 trigonometric polynomial of degree less or equal $ n: $
 $$
Q(n)  = \{ a_0/2 + \sum_{k=1}^n (a_k \cos(kx) + b_k \sin(kx)) \}.
 $$

 \vspace{3mm}
 {\bf D. Auxiliary facts.}\par

 It follows from the classical {\it approximation theory}
 (generalized Jackson inequalities)
 that there exists an absolute constant $ C \in (0,\infty) $ such that

 $$
 E_n(f)_p \le C \cdot  \omega(f,\delta)_p, \eqno(4)
 $$
see for example, \cite{Dai1},  \cite{Ditzan1} .\par

We recall also the well-known Nikol'skii inequality: for all trigonometric polynomial
$ U = U(x) $ of degree less or equal $ n:\ \deg(U) \le n $ and the values
$ p,q: 1 \le p \le q \le \infty $

$$
|U|_q \le 2 \ n^{1/p - 1/q} \ |U|_p, \eqno(5)
$$
see \cite{Ditzan1},  \cite{Nessel1}, \cite{Nessel2}, \cite{Nikol'skii1}.   \par
Note that the inequality (5) is trivially satisfied also in the case $1 \le q < p. $ \par

\section{Bilateral Grand Lebesgue Spaces}

\vspace{3mm}

  We recall briefly the definition and needed properties of these spaces.
  More details see in the works \cite{Fiorenza1}, \cite{Fiorenza2}, \cite{Ivaniec1},
   \cite{Ivaniec2}, \cite{Ostrovsky1}, \cite{Ostrovsky2}, \cite{Kozatchenko1},
  \cite{Jawerth1}, \cite{Karadzov1} etc. More about rearrangement invariant spaces
  see in the monographs \cite{Bennet1}, \cite{Krein1}. \par

\vspace{3mm}

For $a$ and $b$ constants, $1 \le a < b \le \infty,$ let $\psi =
\psi(p),$ $p \in (a,b),$ be a continuous positive
function such that there exists a limits (finite or not)
$ \psi(a + 0)$ and $\psi(b-0),$  with conditions $ \inf_{p \in (a,b)} > 0 $ and
 $\min\{\psi(a+0), \psi(b-0)\}> 0.$  We will denote the set of all these functions
 as $ \Psi(a,b). $ \par
 The Bilateral Grand Lebesgue Space (in notation BGLS) $  G(\psi; a,b) =
 G(\psi) $ is the space of all measurable
functions $ \ f: X \to R \ $ endowed with the norm

$$
||f||G(\psi) \stackrel{def}{=}\sup_{p \in (a,b)}
\left[ \frac{ |f|_p}{\psi(p)} \right], \eqno(6)
$$
if it is finite.\par
 In the article \cite{Ostrovsky2} there are many examples of these spaces.
 For instance, in the case when  $ 1 \le a < b < \infty, \beta, \gamma \ge 0 $ and

 $$
 \psi(p) = \psi(a,b; \beta, \gamma; p) = (p - a)^{-\beta} (b - p)^{-\gamma};
 $$
we will denote
the correspondent $ G(\psi) $ space by  $ G(a,b; \beta, \gamma);  $ it
is not trivial, non - reflexive, non - separable
etc.  In the case $ b = \infty $ we need to take $ \gamma < 0 $ and define

$$
\psi(p) = \psi(a,b; \beta, \gamma; p) = (p - a)^{-\beta}, p \in (a, h);
$$

$$
\psi(p) = \psi(a,b; \beta, \gamma; p) = p^{- \gamma} = p^{- |\gamma|}, \ p \ge h,
$$
where the value $ h $ is the unique  solution of a continuity equation

$$
(h - a)^{- \beta} = h^{ - \gamma }
$$
in the set  $ h \in (a, \infty). $ \par

We will denote for simplicity $ \Psi(1,b) = \Psi(b), \ 1 < b \le \infty; \
\psi(1,b) = \psi(b); \ \psi(1,b;\alpha,\beta) = \psi(b,\beta); $ the value $ \alpha $
in this case is not essential. \par

 The  $ G(\psi) $ spaces over some measurable space $ (X, F, \mu) $
with condition $ \mu(X) = 1 $  (probabilistic case)
appeared in an article \cite{Kozatchenko1}.\par
 The BGLS spaces are rearrangement invariant spaces and moreover interpolation spaces
between the spaces $ L_1(R^d) $ and $ L_{\infty}(R^d) $ under real interpolation
method \cite{Carro1}, \cite{Jawerth1}. \par
It was proved also that in this case each $ G(\psi) $ space coincides
with the so - called {\it exponential Orlicz space,} up to norm equivalence. In others
quoted publications were investigated, for instance,
 their associate spaces, fundamental functions
$\phi(G(\psi; a,b);\delta),$ Fourier and singular operators,
conditions for convergence and compactness, reflexivity and
separability, martingales in these spaces, etc.\par
 Let $ g: X \to R $ be some measurable function such that

 $$
 g \in \cup_{p > 1} L(p).
 $$
We can then introduce the non-trivial function $ \psi_g(p) $ as follows:

$$
\psi_g(p) \stackrel{def}{=} |g|_p.\eqno(7)
$$
This choosing of the function $ \psi_g(\cdot) $ will be called {\it natural
choosing. }\par

{\bf Remark 1.} If we introduce the {\it discontinuous} function

$$
\psi_r(p) = 1, \ p = r; \psi_r(p) = \infty, \ p \ne r, \ p,r \in (a,b) \eqno(8)
$$
and define formally  $ C/\infty = 0, \ C = \const \in R^1, $ then  the norm
in the space $ G(\psi_r) $ coincides with the $ L_r $ norm:

$$
||f||G(\psi_r) = |f|_r.
$$

Thus, the Bilateral Grand Lebesgue spaces are the direct generalization of the
classical exponential Orlicz's spaces and Lebesgue spaces $ L_r. $ \par

The BGLS norm estimates, in particular, Orlicz norm estimates for
measurable functions, e.g., for random variables are used in PDE
\cite{Fiorenza1}, \cite{Ivaniec1}, theory of probability in Banach spaces
\cite{Ledoux1}, \cite{Kozatchenko1},
\cite{Ostrovsky1}, in the modern non-parametrical statistics, for
example, in the so-called regression problem \cite{Ostrovsky1}.\par

\vspace{3mm}

 We use symbols $C(X,Y),$ $C(p,q;\psi),$ etc., to denote positive
constants along with parameters they depend on, or at least
dependence on which is essential in our study. To distinguish
between two different constants depending on the same parameters
we will additionally enumerate them, like $C_1(X,Y)$ and
$C_2(X,Y).$ The relation $ g(\cdot) \asymp h(\cdot), \ p \in (A,B), $
where $ g = g(p), \ h = h(p), \ g,h: (A,B) \to R_+, $
denotes as usually

$$
0< \inf_{p\in (A,B)} h(p)/g(p) \le \sup_{p \in(A,B)}h(p)/g(p)<\infty.
$$
The symbol $ \sim $ will denote usual equivalence in the limit
sense.\par
 The particular orderings in the set $ G\Psi $ may be introduced as follows.
 We will write $ \psi_1(\cdot) < \psi_2(\cdot), $ or equally
 $ \psi_2(\cdot) > \psi_1(\cdot), $ if

$$
\sup_{p \in(A,B)}\psi_1(p)/\psi_2(p)<\infty. \eqno(9)
$$

 We will write $ \psi_1(\cdot) << \psi_2(\cdot), $ or equally
 $ \psi_2(\cdot) >> \psi_1(\cdot), $ if

$$
\lim_{\psi_2(p)\to \infty} \psi_1(p)/\psi_2(p)=0. \eqno(10)
$$

We will denote as ordinary the indicator function
$$
I(x \in A) = 1, x \in A, \ I(x \in A) = 0, x \notin A;
$$
here $ A $ is a measurable set.\par

\bigskip

\section{Main result: upper estimations for embedding theorem}

\vspace{3mm}
 Note at fist that if $ \psi_1(\cdot) < \psi_2(\cdot), $ then $ G(\psi_1) $ is closed
 subspace of the space $ G(\psi_2): $ if $ f \in G(\psi_1) $ and $ ||f||G(\psi_1 = 1,$
 then
$$
|f|_p \le \psi_1(p) \le \psi_2(p)\times \sup_p \frac{\psi_1(p)}{\psi_2(p)} \le C
\times \psi_2(p),
$$
hence $ ||f||G(\psi_2) \le C. $\par
 As a consequence: if $ \psi_1(\cdot) < \psi_2(\cdot) $ and $ \psi_2(\cdot) < \psi_1(\cdot), $ then the spaces $ G(\psi_1) $ and $ G(\psi_2) $ coincide up to
 norm equivalence.\par
 It is known, see \cite{Ostrovsky2}, that if $ \psi_1(\cdot) << \psi_2(\cdot), $
 then the space $ G(\psi_1) $ is compact embedded in the space $ G(\psi_2). $ \par

Further, we suppose $ f(\cdot) \in G(\psi) $ for some $ \psi(\cdot) \in \Psi(b). $
Let us introduce the set $ Z $ as a set of all strict monotonically increasing
sequences of natural numbers $ Z = \{ n(k) \} $ such that
$$
n(1) = 1, \ n(k+1) \ge n(k)  + 1, \ k = 1,2,\ldots. \eqno(11)
$$
For instance, the sequence $ n_0(k)= 2^k - 1 $ belongs to the set $ Z. $  \par
 Let us define a  new function (more exactly, a functional)
 $$
 \nu(q,p) = \nu(||f||G(\psi), \omega^{(G(\psi))}(f,\cdot);\psi; q,p)
 $$
as follows:

$$
\nu(q,p) = \inf_{\{n(\cdot) \in Z \}}
\left\{||f||G(\psi) + 32 \sum_{k=1}^{\infty} \left[ n^{1/p - 1/q}(k+1)
\cdot \omega^{(G(\psi))}(f,1/n(k)) \right] \right\},
$$

$$
\theta(q)= \theta(q;||f||G(\psi), \omega^{(G(\psi))}(f,\cdot);\psi; q,p) =
\inf_{p \in (1,b)}\nu(||f||G(\psi), \omega^{(G(\psi))}(f,\cdot);\psi; q,p). \eqno(12)
$$

{\bf Theorem 1.} {\sc Let} $ \exists \psi(\cdot) \in \Psi(b) $ {\sc for which}
$ f(\cdot) \in G(\psi). $ {\sc Assume that for some} $ b_1 \ge b \  $ {\sc the function}
$ \nu = \nu(q), \ q \in (1,b_1) $ {\sc is finite. Then}

$$
|f|_q \le 6 \theta(q)
$$
{\sc or equally}

$$
||f||G(\theta) \le 6. \eqno(13)
$$
{\bf Proof.} We follow for the method offered by Il'in, see for example \cite{Besov2},
chapter 4.\par
 Let the function $ f = f(x), \ x \in X $ belongs to the space $ G(\psi); $
so that therefore
$$
|f|_p \le \psi(p), \ \sup_{h: |h| \le \delta} |T(h)f - f|_p \le \omega_p(f,\delta) \le
\psi(p) \ \omega^{(G(\psi))}(\delta), \ p \in (1,b).
$$
Let $ n(k) $ be arbitrary sequence from the set $ Z $ and
 let also $ m = m(n) = \Ent(n/2), \ n \ge 4, $ where $ \Ent(z) $ denotes the integer part
 of a positive number $ z, $
 $ V_n[f](x) = V_n(x), n = 2,3,4, \ldots $ be a sequence of the classical
 trigonometrical Valee-Poussin polynomials of degree less than $ 2n $ for the function $ f(\cdot). $ Recall that $  V_n[f](x) = $

 $$
 [2 \pi m(n)]^{-1}\int_0^{2 \pi} f(x+t) \ \sin[(2n+1-m(n))/2] \
 \sin[0.5(m(n)+1)t] \ \sin^{-2}(t/2) \ dt. \eqno(14)
 $$

 It follows from the generalized Jackson inequalities \cite{Achieser1}, \cite{Efimov1},
\cite{Natanson1}, volume 1, chapter 4, p. 83-85, \cite{Timan1} that
$$
|V_n[f]|_p \le C_0 |f|_p,
$$

$$
|f - V_n[f]|_p \le C \omega(f,1/n)_p \le C \ \psi(p) \ \omega^{(G(\psi))}(f, 1/n).
\eqno(15)
$$
with some {\it absolute} constants $C_0, C; $  for example, it may be fetched
 $C_0 = 3, \ C = 12, $ see \cite{Natanson1}, volume 1, chapter 4, p. 83-85.  \par
 Let us consider the following sequence of trigonometrical polynomials

 $$
P(x) = Q_{n(1)}(x); \ Q_k(x) = V_{n(k+1)}[f](x) - V_{n(k)}[f](x); \eqno(16)
 $$
then  $ |P(\cdot)| \le 3 \psi(p) ||f||G(\psi) $ and
$ \deg{Q_k} \le 2 \ n(k+1); $

$$
 |Q_k|_p \le |f -  V_{n(k+1)}[f]|_p + |f -  V_{n(k)}[f]|_p \le
 C \ [\omega(f,1/n(k+1))_p + \omega(f,1/n(k))_p ] \le
$$

$$
2C \ \omega(f,1/n(k))_p \le C_2 \ \psi(p) \ \omega^{(G(\psi))}(f, 1/n(k)),
\eqno(17)
$$
and
$$
f(x) = P(x) + \sum_{k=1}^{\infty} Q_k(x)
$$
almost everywhere and in all the $ L(p) $ sense for all the values $ p $ for
which $ \psi(p) < \infty, $ i.e.
 for all the values $ p $ from the interval $ p \in (1,b). $ \par
We conclude using the Nikol'skii inequality:

$$
|P|_q \le C_2 \psi(p) ||f||G(\psi),
$$

$$
|Q_k|_q \le C_3 \ [n(k+1)]^{1/p - 1/q} \ \psi(p) \ \omega^{(G(\psi))}(f, 1/n(k)).
$$
 We get using the triangle inequality for the $ L(q) $ norm:

$$
|f|_q \le  C_2 \psi(p) ||f||G(\psi) +
C_3 \sum_{k=1}^{\infty}
[n(k+1)]^{1/p - 1/q} \ \psi(p) \ \omega^{(G(\psi))}(f, 1/n(k)). \eqno(18)
$$
 Since the value $ p $ and the sequence $  \{n(k)\} $ are arbitrary, we obtain
 tacking the infimum over $ p $ and $ \{n(k)\} $  the  assertion of theorem 1. \par

\bigskip

\section{Examples.}

\vspace{3mm}
{\bf 1.} If we choose $ n(k) = 2^k - 1, $ we state: $ f(\cdot) \in G(\zeta), $ where

$$
\zeta(q) = \inf_{p \in (1,b)} \left[ \psi(p) +  \sum_{k=1}^{\infty}
2^{k(1/p - 1/q)} \ \psi(p) \ \omega^{(G(\psi))}(f, 2^{-k}) \right].
$$
\vspace{3mm}
{\bf 2.} Further, suppose in addition to the conditions for last  assertion
that there is a constant $ \alpha \in (0, 1/b) $ such that
$$
\omega^{(G(\psi))}(f, \delta) \le C_4 \ \delta^{\alpha}, \ \delta \in [0,1];
$$
then
$$
\zeta(q) = \inf_{p \in (1,b)} \frac{\psi(p) }{\alpha - (1/p - 1/q)},\eqno(19)
$$
as long as

$$
\sum_{k=1}^{\infty} 2^{-k(\alpha -(1/p - 1/q))} \sim (\alpha - (1/p - 1/q))^{-1},
$$
when $ \alpha \in ( 1/p - 1/q, 1/p - 1/q +1), \ \alpha < 1/b. $ \par

Notice that the function $\zeta(q) $ is finite for all the values $ q $ from the
open interval $ q \in (1,b_1), \ b_1 \stackrel{def}{=} b/(1 - \alpha b)). $ \par

{\bf Remark 2.} The assertion (19) is exact, for instance, for the $ G(\psi_r) $
spaces.\par
\vspace{3mm}

{\bf 3.} Furthermore, if  $ \exists \alpha = \const  \in (0, 1/b), \ \alpha_2 =
\const > -1 $ such that
$$
\omega^{(G(\psi))}(f, \delta) \le C_5 \ \delta^{\alpha} \ |\log \delta|^{\alpha_2},
 \ \delta \in [0,1/e],
$$
then
$$
\zeta(q) = \inf_{p \in (1,b)} \frac{\psi(p) }{[\alpha - (1/p - 1/q)]^{\alpha_2 + 1}}.
$$
 When $ \alpha_2 = -1, $ then
$$
\zeta(q) = \inf_{p \in (1,b)}
\{\psi(p) \cdot [ |\log(\alpha - (1/p - 1/q))| + 1] \}.
$$
Note that in the case $ \alpha_2 < - 1 $ we obtain a trivial result:
$ \nu(q) \asymp \psi(q). $ \par
 The last assertions follow from the following elementary fact.\par
 {\bf Lemma 1.} Let the value $ \eta $ belongs to the interval $ (0,1). $ Consider the
 series
 $$
 W(\eta) = W_{\alpha_2}(\eta) = \sum_{k=1}^{\infty} 2^{-k \eta} k^{\alpha_2},
 \alpha_2 = \const.
 $$
 The following equalities are true:

 $$
 {\bf A.} \alpha_2 > -1 \ \Rightarrow \ W_{\alpha_2}(\eta) \asymp \eta^{-1-\alpha_2};
 \eqno(20A)
 $$

$$
 {\bf B.} \alpha_2 = -1 \ \Rightarrow \ W_{\alpha_2}(\eta) \asymp
[ |\log(\eta)| + 1]; \eqno(20B)
$$

$$
 {\bf C.} \alpha_2 < -1 \ \Rightarrow \
 \sup_{\eta \in (0,1)} W_{\alpha_2}(\eta) < \infty. \eqno(20C)
$$
 Note that the last case is trivial for us.\par

\vspace{3mm}

{\bf 4.} Suppose in addition $ \exists b \in (1,\infty), \beta \ge 0, $ such that

$$
\psi(p) = (b - p)^{-\beta},  \ p \in (1,b). \eqno(21)
$$
We will distinguish a two cases: $  \alpha < 1/b $ and $ \alpha  1/b. $ \par
{\bf A1.} The possibility $ \alpha < 1/b $ and $ \alpha_2 < - 1.$ \par
Let us denote $ b_1 = b/(1 - \alpha b). $ \par
 We have after simple calculations for the values $ q \in (1, b/(1-\alpha b)) = (1,b_1) $

$$
\nu(q) \le \zeta_{\beta}(q) \stackrel{def}{=} \inf_{p \in (1,b)}
\frac{(b-p)^{-\beta} }{[\alpha - (1/p - 1/q)]^{\alpha_2 + 1}} \sim
$$

$$
(b_1 - q)^{ -\beta - \alpha_2 - 1}, \ q \in (1, b_1).
$$

 Thus, if $ f(\cdot) \in G\Psi(b; \beta) $ and
 $$
 \omega^{(G\Psi(b;\beta))}(\delta) \le C \ \delta^{\alpha} |\log \delta|^{\alpha_2}, \
 \delta \in (0, 1/e),
 $$
then $ f \in G\Psi(b_1, \beta + \alpha_2 + 1) $ and

$$
||f||G\Psi(b_1, \beta + \alpha_2 + 1) \le C(\alpha,\alpha_2,\beta) \
||f||G\Psi(b, \beta). \eqno(22)
$$

{\bf A2.} The variant  $ \alpha < 1/b, $ but $ \alpha_2 = - 1.$ \par
 Here

$$
\nu(q) \le C \zeta_{\beta,1}(q) \stackrel{def}{=} \inf_{p \in (1,b)}
[(b-p)^{-\beta} ] \cdot[ |\log(\alpha - (1/p - 1/q))| + 1] \sim
$$

$$
(b_1 - q)^{ -\beta } \cdot [|\log(b_1 - q)| + 1 ], \ q \in (1, b_1).
$$

{\bf B1.} Let now $ \alpha = 1/b $ and $ \alpha_2 > -1. $ \par
 Note that in the case if $ \alpha = 1/b $ we conclude formally $ b_1 = + \infty. $ \par
We have using theorem 1:

$$
\nu(q) \le C(\alpha,\alpha_2,\beta) \ q^{\beta + \alpha_2 + 1}, \  q \in (1,\infty).
$$

{\bf B2.} Finally, let $ \alpha = 1/b, \ \alpha_2 = - 1. $ We obtain:

$$
\nu(q) \le C_1(\alpha,\beta) \ q^{\beta} \cdot [\log q + 1].
$$

\bigskip

\section{Exactness of the upper bounds}

\vspace{3mm}

In this section we built some examples in order to illustrate the
 exactness of upper estimations. \par
 {\bf Theorem 2.} {\sc The estimation of theorem 1 is in general case  the best
 possible.}\par
 In order to prove this assertion, we  consider the following example. Let
 $ X = (0,1), \ d \mu = dx, \ \Delta = \const \in (0,1), \
 \gamma = \const > 1, \ b = \const \in (1,1/\Delta);$  we will suppose further
 $ \gamma \to \infty; $

 $$
 f_0(x):= x^{-\Delta} \ |\log x|^{\gamma} \ I(x \in (0,1)). \eqno(23)
 $$

We find by the direct calculations at $ p \in (1,1/\Delta):$

$$
|f_0|_p^p = \int_0^1 x^{-p \Delta} |\log x|^{\gamma p} d x =
\frac{\Gamma(\gamma p + 1)}{(1 - p \Delta)^{\gamma p + 1} }, \eqno(24)
$$

$$
|f_0|_p \asymp (q_0 - p)^{-\gamma - 1/p}, \ q_0 \stackrel{def}{=} 1/\Delta.
$$
and we take the natural function
$$
\psi(p) = |f_0|_p, \ p \in (1,b);
$$
then

$$
\omega^{(G(\psi_0))}(f_0,\delta) \asymp \delta^{1/b - 1/q_0} \ |\log \delta|^{\gamma}, \
\delta \in (0, 1/e); \eqno(25)
$$
so that here

$$
\alpha = 1/b - 1/q_0, \ \alpha_2 = \gamma, \ \beta = 0.
$$
It follows from theorem 1 that $ f_0 \in G\nu(q_0,1/\Delta), $ where
$$
\nu(q) \asymp (q_0 - q)^{-\gamma - 1}, \ q \in (1, q_0) = (1, 1/\Delta);
$$
but really

$$
\nu_0(q) \stackrel{def}{=} |f_0|_q \asymp (q_0 - q)^{-\gamma - \Delta}.
$$
We can see that the {\it extremal value} of the function $ \nu_0, $ i.e. the value
$$
q_{\max} \stackrel{def}{=} \sup \{ q, \ \nu_0(q) < \infty  \} \eqno(26)
$$
is calculated exactly by means of theorem 1: $ q_{\max} = q_0 = 1/\Delta. $ \par
 The {\it power} of the function $ \nu_0, $ i.e. the value

 $$
 w(\nu_0) \stackrel{def}{=} \overline{\lim}_{q \to q_0-0}
 \frac{|\log \nu_0(q)|}{|\log (q - q_0)|} \eqno(27)
 $$
evaluating by means of theorem 1 is equal to $ \gamma + 1, $ but really this value
is equal to $ \gamma + \Delta. $ \par
 Notice that

 $$
 \lim_{\gamma \to \infty} \frac{\gamma + \Delta}{\gamma + 1} = 1. \eqno(28)
 $$

\vspace{3mm}

\section{Concluding remark: compactness or not compactness of embedding operator}

\vspace{2mm}
Let $ \psi(\cdot) $ and $ \theta(\cdot) $ be at the same as in theorem 1. Note that 
in general case the unit embedding operator $ T: G(\psi) \in G(\theta) $ is not compact
\cite{Ostrovsky2}. \par
 But if we consider the unit operator $ T: G(\psi) \in G(\xi), $ where $ \xi(\cdot) 
 << \psi(\cdot), $ then this operator $ T $ is compact.\par

\end{document}